\documentclass[reqno]{amsart}
\usepackage{graphicx}
\usepackage[all]{xy}
\theoremstyle{theorem}
\newtheorem{theorem}{Theorem}[section]

\newtheorem{lemma}[theorem]{Lemma}
\newtheorem{prop}[theorem]{Proposition}

\theoremstyle{definition}
\newtheorem{definition}{Definition}[section]
\newtheorem{remark}[definition]{Remark}

\numberwithin{equation}{section}

\newcommand{\CC}{\mathbb{C}}
\newcommand{\RR}{\mathbb{R}}
\newcommand{\QQ}{\mathbb{Q}}
\newcommand{\kk}{\boldsymbol{k}}
\newcommand{\ZZ}{\mathbb{Z}}

\newcommand{\AI}{A_\infty}
\newcommand{\LI}{L_\infty}
\newcommand{\CP}{\mathbb{C}P}

\newcommand{\CM}{\mathcal{M}}

\newcommand{\E}{\epsilon}
\newcommand{\HH}[1]{\widehat{#1}}
\newcommand{\WH}[1]{\widehat{#1}}
\newcommand{\WT}[1]{\widetilde{#1}}
\newcommand{\OL}[1]{\overline{#1}}
\newcommand{\UL}[1]{\underline{#1}}
\newcommand{\NOV}{\Lambda_{nov}}
\newcommand{\NOVO}{\Lambda_{0,nov}}
\begin{document}
\title{On the counting of holomorphic discs \\ in toric Fano manifolds}

\author{Cheol-Hyun Cho}
\address{Cheol-Hyun Cho, Department of Mathematics and Research Institute of Mathematics, Seoul National University,
Kwanakgu, Seoul, South Korea, Email: chocheol@snu.ac.kr}

\begin{abstract}
Open Gromov-Witten invariants in general are not well-defined. We discuss in detail the enumerative numbers of the Clifford torus $T^2$ in $\CP^2$.
For cyclic A-infinity algebras, we show that certain generalized way of counting may be defined up to Hochschild or cyclic boundary elements. In particular
we obtain a well-defined function on Hochschild or cyclic homology of a cyclic A-infinity algebra, which
has invariance property under cyclic A-infinity homomorphism. We discuss an example of Clifford torus $T^2$ and compute the invariant for
a specific cyclic cohomology class.
\end{abstract}
\thanks{This work was supported by the SRC Program of Korea Science and Engineering Foundation (KOSEF) grant funded by the Korean government(R11-2007-035-01001-0)}
\maketitle
\bigskip
\section{Introduction}
Gromov-Witten invariants, which are obtained from the intersection theory of the moduli spaces of pseudo-holomorphic curves without boundary, have provided exciting links to various branches of mathematics. Open Gromov-Witten invariants
 from the intersection theory of pseudo-holomorphic curves with boundary, are still largely mysterious. 
Most of the works in this direction has been carried out in the case that 
Lagrangian submanifold is given as a fixed point set of anti-symplectic involution or in the Calabi-Yau case
(See for example, Katz-Liu \cite{KL}, Welschinger \cite{W}, Cho \cite{C4}, Solomon \cite{S},
 Pandharipande-Solomon-Walcher \cite{PSW}, Fukaya \cite{Fu2})

The fundamental difference between open and closed cases is that while the moduli spaces of closed pseudo-holomorphic curves carry fundamental cycles, but the moduli spaces of open pseudo-holomorphic curves do not carry fundamental cycles but only fundamental chains. As the chains do not have good intersection property, it is not clear in general how to define open Gromov-Witten invariants. For example, even in the case of anti-symplectic involution, it is not yet known how to define them as an invariant for dimension $> 3$. 

Even though they are not invariants, there has been several approaches to obtain certain invariant from it as
in the work of Biran-Cornea, and more recently by Fukaya-Oh-Ohta-Ono \cite{FOOO}. 
In this paper, we show two approaches: one by considering a three point intersection problem where we take a combination of chains to make a Floer cycle. Another more general approach is to take a combination of intersection problems to make it a cycle of Hochschild or cyclic homology.

In the latter approach, we try to find a combination of intersection problems whose total sum becomes invariant. 
For this purpose, we investigate the condition on the sum of chains to be intersected at the boundary
so that totality of the bubbling off a disc becomes zero. 
The simple idea is that this condition is related to Hochschild or cyclic boundary operations.

Using the formalism of cyclic $\AI$-algebras, we
 show that for a cyclic A-infinity algebras, such counting may be defined up to Hochschild or cyclic boundary elements, which may provide an invariant counting in a generalized sense. In particular 
we obtain a well-defined function on Hochschild or cyclic homology of a cyclic A-infinity algebra, which
has invariance property under cyclic A-infinity homomorphism. 

In the last section, we discuss the example of Clifford torus in $\CP^2$, where we find a specific cyclic cohomology cycle of the
$\AI$-algebra of the Clifford torus $T^2$ and compute the invariant defined in this paper.

This is a revised version of the unpublished manuscript arXiv:math/0604502 using the language of cyclic $\AI$-algebras. We thank the
anonymous referee for pointing out a much shorter (and better) proof of Lemma 2.1 and helpful comments to improve the paper.

\section{Example of three point open Gromov-Witten invariant in Clifford tori}
We examine the enumerative numbers of holomorphic discs with boundary on torus fibers of toric manifolds. It  is meaningful to examine the case of Clifford tori in projective spaces to understand the general feature of enumerative numbers in the open case. 

In particular, in this section, we examine the case of $T^2$ in $\CP^2$, the number of holomorphic discs intersecting three generic points, and also the number of such intersecting discs which preserves the prescribed cyclic order. The latter is related to $m_2$ operation of the $\AI$-algebra and we observe a relation to the work of Biran-Cornea.

Let $M$ be a toric Fano manifold. The Floer theory of Lagrangian torus fibers of $M$ has been
developed in \cite{C1},\cite{CO} and \cite{FOOO1}. The main characteristic of the theory is that its homology and the product structure is determined by the holomorphic discs of Maslov index two.
But the number of holomorphic discs intersecting generic three points has not been discussed there
as it involves Maslov index 4 discs.

We first consider the case of $(\CC^n,(S^1)^n)$, which
can be extended to other toric cases in a standard way.
For a Lagrangian submanifold $L \subset M$, we
denote by $\CM_3(\beta,J)^*$ the moduli space of simple $J$-holomorphic discs
with three boundary marked points of homotopy class $\beta \in \pi_2(M,L)$.
Recall that the moduli space $\CM_3(\beta,J)^*$  has an
expected dimension 
$$dim(L) + \mu(\beta) -dim(PSL(2;\RR)) + 3 = n + \mu(\beta).$$
Here $\mu(\beta)$ is the Maslov index. 
We have an evaluation map $ev_3:\CM_3(\beta,J)^* \to L^{3}$ given by
the evaluation at the marked points. We want to count the
number of points in $ev_3^{-1}(p_0,p_1,p_2)$ where 
$p_i$'s are disjoint generic points in $L$. Hence, to have a zero dimensional
preimage, we require that $\mu(\beta) = 2n$.

Consider $\CC^n$ and a torus $L$ defined by
$$L := \{ (z_1,\cdots,z_n)| \; \forall i , \, |z_i| = 1 \}.$$
Note that $\pi_2(\CC^n, L)$ is generated by $n$
elements of whose Maslov indices are two, which we denote by $\beta_1,\cdots,\beta_n$. They
are homotopy classes of the maps $w_i :D^2 \to \CC^n$ defined
by $z \mapsto (1,\cdots,1,z,1,\cdots,1)$ where $z$ is located
at the $i$-st entry. 

Let us fix a homotopy class of holomorphic disc to be $\sum_{j=1}^n\beta_{j}$,
It is easy to see that for other homotopy classes with $\mu =2n$,
there wouldn't be any intersection between holomorphic discs of such class and generic three points
as one of the factor in $\CC^n$ above has to remain constant.

The holomorphic disc of class  $\sum_{j=1}^n\beta_{j}$ may be written as 
$$(e^{c_1' i}\frac{z-\alpha_1'}{1-\overline{\alpha_1'}z},\cdots,
e^{c_n' i}\frac{z-\alpha_n'}{1-\overline{\alpha_n'}z}).$$

We observe that
\begin{lemma}\label{pr:3pt}
Let $$p = (e^{\theta_{p1}i},\cdots,e^{\theta_{pn}i}),
q = (e^{\theta_{q1}i},\cdots,e^{\theta_{qn}i}),
r = (e^{\theta_{r1}i},\cdots,e^{\theta_{rn}i})$$
be generic three points on the torus $L \subset \CC^n$. 
Then if the cyclic ordering of three complex numbers  
$(e^{\theta_{pj}i},e^{\theta_{qj}i},e^{\theta_{rj}i})$ on $S^1$ for each $j$ are
all the same (i.e. all counter-clockwise or all clockwise), then there is a unique holomorphic disc 
$(D^2,\partial D^2) \to (\CC^n,L)$ of the homotopy class 
$\sum_{j=1}^n \beta_j$, which passes through $p, q, r$.
Otherwise, there does not exist a holomorphic disc of the given homotopy class which passes through  these three points.
\end{lemma}
\begin{proof}
In the preprint version, we have given a proof based on the plane geometry. The following
proof, which was suggested by an anonymous referee, is much shorter and easier. 
Consider a disc $D^2 \subset \CC$ with three fixed marked points $1,i,-1$.
Then, it is a standard fact that $PSL_2(\RR)$ act transitively on triples of cyclically ordered points in the boundary of the disc.  Hence, if 
$(e^{\theta_{pj}i},e^{\theta_{qj}i},e^{\theta_{rj}i})$ on $S^1$ for each $j$ are all counter-clockwise then
we can find a holomorphic map $u_j:D^2 \to D^2$ mapping $1,i,-1$ to the above three points for each $j$.
Then, the desired holomorphic disc passing through $p, q, r$ is nothing but
$u= (u_1,\cdots,u_n)$.
In the clockwise case, consider the ordered triple $p, r, q$ instead to obtain a similar result.
In the case the cyclic ordering for each $j$ is not the same for all $j$, it is easy to see that existence
of a holomorphic disc would violate the standard fact that $PSL_2(\RR)$ preserves the cyclic orderings.
\end{proof}

In a joint work \cite{CO} with Yong-Geun Oh, we have classified all holomorphic discs with 
boundary on Lagrangian torus fibers in toric manifolds and showed that the discs in the homogeneous coordinate, are
given by standard Blaschke products. Hence the above result can be used if the intersection occurs
in one of the standard open set $\CC^n$ of the toric variety. We also remark that 
 the standard complex structure $J_0$ is Fredholm regular (as shown in \cite{CO}) and hence this may be considered as a generic phenomenon.

Before we state the main result of this section,
we first explain a two dimensional chain $Q_p^J$ which depends on $p$ and $J$.
In \cite{C2}, we explained that although Bott-Morse Floer cohomology of the Clifford torus is isomorphic to the singular cohomology of the torus, the actual cycles of each homology are different. For example $\partial p =0$ in the singular cohomology for a point $p=[p_0;p_1;p_2] \in T^2$, but
$$\delta_{HF}(p) = l_0 + l_1+l_2,$$
where $l_0,l_1,l_2$ are the boundary images of  the holomorphic discs $[p_0z;p_1;p_2]$, $[p_0;p_1z;p_2]$, $[p_0;p_1;p_2z]$.
Note that $l_0 +l_1+l_2$ is not zero on the chain level,
but is homologous to zero. Hence we can choose a 2-chain $Q_p^{J_0}\subset L$ with 
$ \partial Q_p^{J_0} =-( l_0 +l_1+l_2)$, then  the sum
$p + Q_p^{J_0}\otimes T^{\omega(D)}$, turns out to be the correct cycle in Floer homology (see \cite{C2} for details). We set $\omega(D)=1$ for simplicity from now on.

Now, here is the main result in this section
\begin{prop} Let $L$ be Clifford torus in $\CP^2$. In what follows, by disc, we mean a holomorphic disc in $\CP^2$ with boundary on $L$. Let $p, q, r$ be three generic points on $L$. 
\begin{enumerate}
\item The number of  discs passing through $p, q, r$ of a {\em fixed homotopy class} depends on the choice of $p, q, r$ and hence is not an invariant.
\item The number of  discs  passing through $p, q, r$
is even.
\item The number of  discs  passing through $p, q, r$ 
in a cyclic order depends on the choice of $p, q, r$  and is not an invariant.
\item The number in $(3)$ can be made invariant by adding another intersection number $T_{pqr}$.
\end{enumerate}
\end{prop}
The number $T_{pqr}$ is defined to be the sum of intersection numbers of
$$Q_p^{J_0} \cap Q_q^{J_0} \cap r, \;\;Q_p^{J_0} \cap q \cap Q_r^{J_0}, \;\;  p \cap Q_q^{J_0} \cap Q_r^{J_0}.$$
\begin{proof}
We may set $p=[1;1;1]$ by torus action as all ingredients of the proposition admit torus action.
We first denote the homotopy classes of holomorphic discs of Maslov index two as $\beta_0,\beta_1,\beta_2$
each of which corresponding to $[z;1;1],[1;z;1],[1;1;z]$.

Note that  there are three homotopy classes of Maslov index 4 discs, $\beta_1+\beta_2, \beta_0+\beta_2, \beta_0+\beta_1$
which contributes to the three point intersections.
The part (1) is easy to prove: Fix a homotopy class to be $\beta_1+\beta_2$ without loss of generality. Then holomorphic disc does not intersect the divisor $[0;z_1;z_2]$ (see \cite{CO}) and hence lies inside
the chart $U_0=\{[1;z_1;z_2]\}$. Then the previous lemma proves the part (1). This can be easily seen from the figure. We draw $L=T^2$ as
a square with opposite sides identified  and take $p$ to be the vertex. If we fix $q$ inside the square, the region of $r$ which 
intersects a holomorphic disc passing through $p, q$ can be drawn as a shaded region in Figure 1 (0).

\begin{figure}[h]
\begin{center}
\includegraphics[height=1.3in]{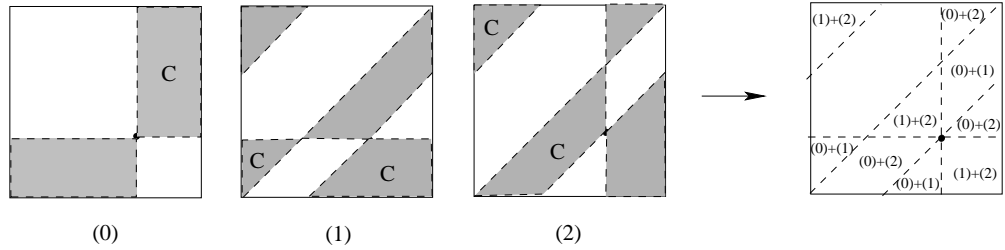}
\caption{The case of Clifford torus in $\CP^2$}
\label{cp2}
\end{center}
\end{figure}

To prove (2), we consider all the above three homotopy classes of Maslov index 4.
In the case of $\beta_0+\beta_2$, disc lies in $U_1=\{[w_1;1;w_2]\}$ with
$w_1 = 1/z_1, w_2 = z_2/z_1$. The cyclic condition of $p, q, r$ in terms of coordinate $w_1, w_2$ is translated to the
condition of $z_1, z_2$ coordinates of $p, q, r$ as follows: namely, using the notation of the previous lemma,
the points $p, q, r$ in given in the coordinates of $U_0 \cong \CC^n$ meets
a disc of homotopy class $\beta_0 + \beta_2$ if the cyclic ordering of 
$1, 1/q_{1}, 1/r_1$ agrees with the cyclic ordering of 
$1,  q_{2}/q_{1},r_{2}/r_{1}$. In the square figure (with coordinate $x,y$),
these orderings are measured by lines $x=c_1$ and $y-x=c_2$. In particular, we consider
whether the resulting $c_1$ and $c_2$ values of $p, q, r$ are cyclic ordered in $\RR/\ZZ$
in the same way. This condition can be checked readily and drawn as a shaded region of Figure 1, (1).
One can proceed similarly to obtain Figure 1, (2). Hence, the total number can be observed by
considering a sum of three figures in the left hand side of Figure 1, which is even as in every region
either there exist two discs or no disc which intersect the three points.
We learned the vanishing of $\ZZ/2\ZZ$ counting from the discussion with Biran and Cornea who considered somewhat different enumerative invariants for $T^2$ with a different methods in \cite{BC}.
Unfortunately, this kind of phenomenon does not seem to occur for intersection problems with more than 3 points.

To prove (3), we observe that the region of $r$ preserving the cyclic order can be also drawn similarly and is given by the shaded region which is labeled by $C$ in the
Figure 1. By adding up the shaded region labeled by $C$, 
 one can easily note that the region  is not the whole $T^2$ nor doubly covered. In fact the number of such discs equals one  except in the empty region of the right hand side of Figure 1.
 
There are two ways to prove (4). 
Let us first describe  an elementary but tedious way: note that the chain $Q_p^{J_0}$ can be chosen
as an upper triangle obtained by bisecting the square along the diagonal ($l_0$).  One can choose triangles $Q_q^{J_0}, Q_r^{J_0}$ similarly. Hence the contribution
$T_{pqr}$ can be checked by hand and as a result we obtain that $T_{pqr}$ equals one except
in the  the empty region of the right hand side of Figure 1. Hence, comparing the result of (3), 
we can see that the mod 2 number of the sum remains invariant. 

  We describe more conceptual approach which explains that the signed sum should remain invariant. The idea is that the sum of discs in (3) and the number $T_{pqr}$ should be considered as an intersection number with three Floer cycles of points.

Here, the appearance of $T_{pqr}$ is deeply related to the fact that in Bott-Morse Lagrangian Floer homology of $L$, the point $p$ is not a cycle but $p + Q_p^{J_0}T$ is a cycle. Hence to obtain
a suitable invariant, one should take intersection with $p + Q_p^{J_0}T$  not with $p$. Hence, this is a generalized intersection problem where we take a combination of chains to make it a Floer cycle. Later in this paper, we make another approach where we take a combination of intersection problems to make it a cycle of Hochschild or cyclic homology.

Consider the intersection problem of holomorphic discs with $p + Q_p^{J_0}T$,
$q + Q_q^{J_0}T$, $r + Q_r^{J_0}T$. This should be suitably interpreted as $p$ is of dimension 0 but
$Q_p$ is of dimension two. Hence, by the intersection number, we mean the sum of following numbers
\begin{enumerate}
\item[A.] The number of holomorphic discs of Maslov index four passing through $p,q,r$ with a given order.
\item[B.] The number of holomorphic discs of Maslov index two passing through two points  out of
$p,q,r$ and corresponding two chain $Q$ of the other point.
\item[C.] The intersection number of 
$$Q_p^{J_0} \cap Q_q^{J_0} \cap r, \;\;Q_p^{J_0} \cap q \cap Q_r^{J_0}, \;\;  p \cap Q_q^{J_0} \cap Q_r^{J_0}.$$
\end{enumerate}
In fact the number of discs in B always vanishes as Maslov index two discs does not pass through
generic two points. Hence, the standard intersection number of three Floer cycles corresponding to
a point equals $A +C$.

Once can check that this number is independent of the choice of $p,q,r$. 
Notice that  the bubbling off
of disc from the Maslov index 4 disc will be canceled by the intersection problem with
$\partial Q$. This is very similar to the technique developed  in the obstruction theory of \cite{FOOO}.
We leave the details to the reader.
This proves $(4)$.
\end{proof}
\begin{remark}
In fact, one can check that the number $T_{pqr}$ is strongly related to
the number $n_qn_r$ (Biran-Cornea \cite{BC2} page 41). Biran and Cornea defined the numbers $n_p, n_q,  n_r$ as
follows. Choose a triangle $\Delta$ with vertex $p, q, r$. Denote the number of holomorphic discs through $p$ intersecting an edge $\overline{qr}$ of a triangle $\Delta$ is denoted as $n_p$. The numbers $n_q, n_r$ are defined similarly.
In this example of Clifford torus $T^2$, by tedious but elementary inspection of the configurations prove that  $T_{pqr} + n_qn_r \equiv 1  (mod 2)$ independent of  any choice of $p, q, r$ or the triangle. 
Hence the invariance of the number (4) in our proposition is equivalent to the invariance result of
the number of discs of (3) together with $n_qn_r$ in \cite{BC2}, which is proved in a different way.
This provides a hint to understand the mysterious counting $n_p,n_q,n_r$ obtained by Biran and Cornea. It would be interesting to find other relations as above to interpret the invariants obtained by them. 
\end{remark}
\begin{remark}
The invariance of (4) can be extended to any other toric Fano case of dimension two, or the case of 
Clifford torus in $\CP^n$. In the case of dimension two, one can choose the chain $Q$ as above and
the result holds true. In the case of $\CP^n$, there is a general formula for higher order correction terms
in \cite{C2}:
\begin{equation}\label{qqqp}
\Psi(P)= P + Q \times P + \frac{1}{2} Q \times Q \times P +
\cdots + \frac{1}{k} Q \times \cdots \times Q \times P,
\end{equation}
where $Q= Q_{p}^{J_0}$ is the 2 dimensional chain 
defined before with $p = [1,\cdots,1] \in \CP^n$,
and $k$ is the largest integer smaller than $(n-dim(P))/2$.
There is a slight error in the definition 4.4 of \cite{C2}
where we need to replace $\sum_{i_1 < \cdots <i_k}$ by $1/k \sum_{\forall i_1,\cdots,i_k}$.
And the chains such as $Q \times P$ are defined as geometric chains using by multiplication
in $T^n \cong (S^1)^n \subset \CC^n$. 
But in general toric Fano case, we do not know how to define an exact Floer cycle.
\end{remark}

\section{Cyclic $\AI$-algebras}
In this section we recall the definition of cyclic $\AI$-algebras and
also describe another formulation of such structure.  For the definition of $\AI$-algebra and
its homomorphisms, we refer readers to \cite{FOOO}. We remark that in \cite{Fu1}, Fukaya has proved that
the $\AI$-algebra of Lagrangian submanifold can be constructed as cyclic filtered $\AI$-algebra.
 For the rest of the paper, we assume that $\AI$-algebras to be considered are finite dimensional.
\subsection{Cyclic $\AI$-algebras}
The cyclic structure on an $\AI$-algebra was first considered by Kontsevich \cite{K} as an invariant symplectic form on the non-commutative formal manifolds. There are somewhat different sign conventions (whether one is working with degree shifting or not) and
we refer readers to \cite{C3} or \cite{CL} for more details. Let $\kk$ be a field of characteristic 0.
\begin{definition}
An $\AI$-algebra $(A,\{m_*\})$ is said to have a {\it cyclic} inner product if
there exists a skew-symmetric non-degenerate, bilinear map $$<,> : A[1] \otimes A[1] \to \kk,$$
such that for all integer $k \geq 1$,
\begin{equation}\label{cyeqn}
    <m_{k,\beta}(x_1,\cdots,x_k),x_{k+1}> = (-1)^{K}<m_{k,\beta}(x_2,\cdots,x_{k+1}),x_{1}>.
\end{equation}
 where $K = |x_1|'(|x_2|' + \cdots +|x_{k+1}|')$. Here $|x|'=|x|-1$ is the shifted degree of $x$.
 For short, we will call such an algebra, cyclic $\AI$-algebra.
\end{definition}
There is a notion of cyclic $\AI$-homomorphism due to Kajiura \cite{Kaj}
\begin{definition}\label{def:homo}
An $\AI$-homomorphism $\{h_k\}_{k\geq 1 }$ between two cyclic $\AI$-algebras is called
a cyclic $\AI$-homomorphism if
\begin{enumerate}
\item $h_1$ preserves inner product $<a,b> = <h_1(a),h_1(b)>$.
\item \begin{equation}
\sum_{i+j=k} <h_i(x_1,\cdots,x_i), h_j(x_{i+1},\cdots,x_k)> =0.
\end{equation}
\end{enumerate}
\end{definition}

\subsection{Other forms of cyclic symmetry}
Now, we prove certain equations induced from $\AI$-equations and cyclicity.
We first mention a version of Stoke's theorem, which follows immediately from the cyclic symmetry.
\begin{lemma}
$$<m_1(x),y> = (-1)^{|x|} <x,m_1(y)>$$
\end{lemma}
\begin{proof}
$$<m_1(x),y> = (-1)^{|x|'|y|'}<m_1(y),x> =
(-1)^{|x|'|y|' + |x|'(|y|'+1) +1}<x,m_1(y)> $$
where the last equality is from the graded skew symmetry of $<,>$. This proves the
lemma.
\end{proof}

Also $\AI$-equation give rise to the following equality.
\begin{prop}\label{prop:cycplus1}
Let $(A,m_*)$ be a cyclic $\AI$-algebra. We have the following identity.
\begin{equation}\label{formula1}
0 = \sum_{\sigma,k_1,k_2}
(-1)^{Kos }
<m_{k_1}(x_{\sigma(1)}, \cdots,x_{\sigma(k_1)}), m_{k_2}(x_{\sigma(k_1+1)}, \cdots,x_{\sigma(k+1)})>
\end{equation}
where the summation is over all cyclic permutations $\sigma \in \ZZ /(k+1)\ZZ$ and for
$ k_1+k_2=k+1$ with additional condition that $1 \in \{\sigma(1),\cdots,\sigma(k_1)\}$.
Here $(-1)^{Kos }$ is the Koszul sign occurring to
rearrange
$$m_{k_1}m_{k_2}x_1x_2\cdots x_{k+1} \Rightarrow m_{k_1}x_{\sigma(1)}\cdots x_{\sigma(k_1)}m_{k_2}x_{\sigma(k_1+1)}\cdots x_{\sigma(k+1)},$$
where we regard  $m_{*}$ as a letter with degree $1$ and $x_i$ with degree $|x_i|'$.
\end{prop}
\begin{remark}
The additional condition, $1 \in \{\sigma(1),\cdots,\sigma(k_1)\}$, make sure that  $x_1$ always appear on the
first factor. We may remove this condition and then the resulting equation is just the twice of the above equation.
\end{remark}
\begin{proof}
The sign $(-1)^{Kos}$ will be used as defined in the above statement.

First, note that in the expression (\ref{formula1}), the cases when either $k_1$ or $k_2$ equals $1$ is
\begin{equation}\label{eq1_4}
(-1)^{|x_1|'}<m_1(x_1),m_k(x_2,\cdots,x_{k+1})>
\end{equation}
\begin{equation}\label{eq1_5}
(-1)^{Kos}<m_k(x_{i+1},\cdots,x_{i-1}),m_1(x_i)>
\end{equation}
for $i=2,\cdots,k+1$.

We start with the Stoke's theorem applied to $m_k(x_1,\cdots,x_k)$:
\begin{equation}\label{eq1}
-< m_1(m_k(x_1,\cdots,x_k)),x_{k+1}>
$$ $$+ (-1)^{|m_k(x_1,\cdots,x_k)|} <m_k(x_1,\cdots,x_k), m_1(x_{k+1})>
=0\end{equation}
Note that $|m_k(x_1,\cdots,x_k)|=|x_1|'+ \cdots + |x_k|'+2$.
Hence, the second term appears in the expression (\ref{formula1}) with
the correct sign, which is in fact the term (\ref{eq1_5}) when $i=k+1$.

Now, we apply $\AI$-equation to the first term of (\ref{eq1}):
$$-< m_1(m_k(x_1,\cdots,x_k)),x_{k+1}> = $$
\begin{equation}\label{eq1_2}
< m_k( m_1(x_1),\cdots,x_k),x_{k+1}>  + \cdots + (-1)^{Kos}  < m_k(x_1,\cdots, m_1(x_k)),x_{k+1}>
\end{equation}
\begin{equation}\label{eq1_3}
+ \sum_{k_1+k_2=k+1} (-1)^{Kos} <m_{k_1}(x_1,\cdots,x_{i},
m_{k_2}(x_{i+1},\cdots,x_{i+k_2}),\cdots,x_{k}),x_{k+1}>
\end{equation}
One may check that (\ref{eq1_2}) gives rise to the terms in (\ref{eq1_4}) and (\ref{eq1_5}) with
the correct sign. Hence we focus on the expression (\ref{eq1_3}).
Note that they can be divided into two cases, $i=0$ or $i\neq 0$.

In the case that $i \neq 0$, by applying the
cyclic symmetry, we move the expression $m_{k_2}(x_{i+1},\cdots,x_{i+k_2})$ to
the right hand side of the $<,>$:
$$\sum_{\stackrel{k_1+k_2=k+1}{\beta= \beta_1 + \beta_2},\beta_i \neq 0} (-1)^{Kos} <m_{k_1}(x_{i+k_2+1},
,\cdots, x_{k+1}, x_1,\cdots, x_{i}),
m_{k_2}(x_{i+1},\cdots,x_{i+k_2})>. $$

For the case that $i=0$, we have
$$\sum (-1)^{Kos} <m_{k_1}(m_{k_2}(x_1,\cdots,x_{k_2}),x_{k_2+1},\cdots,x_k),x_{k+1})>$$
\begin{eqnarray*}
\;&=& \sum (-1)^{\E_1} <m_{k_1}(x_{k_2+1},\cdots,x_k,x_{k+1}), m_{k_2}(x_1,\cdots,x_{k_2})> \\
&=& \sum (-1)^{\E_2} < m_{k_2}(x_1,\cdots,x_{k_2}),m_{k_1}(x_{k_2+1},\cdots,x_k,x_{k+1}))> \\
&=& \sum (-1)^{\E_3} < m_{k_1}(x_1,\cdots,x_{k_1}),m_{k_2}(x_{k_1+1},\cdots,x_k,x_{k+1}))>.
\end{eqnarray*}
Here, we have
\begin{eqnarray*}
\E_1 &=& 1 \cdot (|x_{k_2+1}|' + \cdots + |x_{k+1}|') + (|x_1|'+\cdots + |x_{k_2}|')(|x_{k_2+1}|' + \cdots + |x_{k+1}|') \\
\E_2 &=& \E_1 + 1 +
(|x_1|'+\cdots + |x_{k_2}|'+1)(|x_{k_2+1}|' + \cdots + |x_{k+1}|'+1)\\
&=& |x_1|'+\cdots + |x_{k_2}|' \\
\E_3 &=& Kos
\end{eqnarray*}
Here, the first equality is from cyclic symmetry and the second equality is from the graded skew symmetry of $<,>$ and
the third equality is just relabeling of indices and the resulting sign exactly corresponds to the required the Koszul sign
of the formula (\ref{formula1}). One may check easily that we have produced all the terms in (\ref{formula1}) with the correct sign.

To prove the statement in the remark that removing the condition on $1$ gives rise to the twice of (\ref{formula1}),
note that the case with $x_1$ on the right hand side of $<,>$ equals the similar term with $x_1$ on the
left hand side of $<,>$ by applying the skew symmetry of $<,>$. In this case these two expressions do not
cancel out as the negative sign of skew symmetry of $<,>$ is canceled out by the negative sign of
exchange of two $m$'s.
\end{proof}
Conversely, given the maps $\{m_k\}$ with cyclic symmetry (\ref{cyeqn})(which is non-degenerate), the equation (\ref{formula1})
is equivalent to the $\AI$-equation.

\subsection{Filtered $\AI$-algebras}
We recall the notion of gapped filtered $\AI$-algebra, and we refer readers to \cite{FOOO} for full details.
To consider $\AI$-algebras arising from the study of Lagrangian submanifolds or in general pseudo-holomorphic curves,
one considers filtered $\AI$-algebras over Novikov rings, where the filtration is given by the energy of pseudo-holomorphic curves.
Here Novikov rings are, for a ring $R$
(here $T$ and $e$ are formal parameters)
$$\NOV = \{ \sum_{i=0}^\infty a_i T^{\lambda_i}e^{q_i} | \; a_i \in R,\;\lambda_i \in \RR,\; q_i \in \ZZ, \; \lim_{i \to \infty} \lambda_i = \infty \}$$
$$\NOVO = \{ \sum_i a_i T^{\lambda_i}e^{q_i} \in \NOV | \lambda_i \geq 0 \}.$$

The gapped condition is defined as follows.
The monoid $G \subset \RR_{\geq 0} \times 2 \ZZ$ is assumed to satisfy the following conditions
\begin{enumerate}
\item The projection $\pi_1(G) \subset \RR_{\geq 0}$ is discrete.
\item $G \cap (\{0\} \times 2\ZZ) = \{(0,0)\}$
\item $G \cap (\{\lambda \} \times 2\ZZ)$ is a finite set for any $\lambda$.
\end{enumerate}

Consider a free graded $\NOVO$ module $C$, and let $\OL{C}$ be an $\kk$-vector space
such that $C = \OL{C} \otimes_{\kk} \NOVO$. Then $(C,m_{\geq 0})$ is said to be $G$-gapped if there exists
homomorphisms $m_{k,\beta}:(\overline{C}[1])^{\otimes k} \to \overline{C}[1]$ for $k=0,1,\cdots,$ $\beta=(\lambda(\beta),\mu(\beta)) \in G$
such that
$$m_k = \sum_{\beta \in G} T^{\lambda(\beta)}e^{\mu(\beta)/2} m_{k,\beta}.$$

Recall that these $m_k$ operations may be considered as coderivations by defining
\begin{equation}\label{eq:hatdfil}
\WH{m}_k(x_1 \otimes \cdots \otimes x_n) = \sum_{i=1}^{n-k+1}
(-1)^{|x_1|' + \cdots + |x_{i-1}|'} x_1 \otimes \cdots \otimes m_k(x_i,
\cdots, x_{i+k-1}) \otimes \cdots \otimes x_n
\end{equation}
for $k \leq n$ and $\WH{m}_k(x_1 \otimes \cdots \otimes x_n) =0$ for $k >n$.
If we set $\WH{d} = \sum_{k=0}^\infty \WH{m}_k$, the $\AI$-equations are equivalent to the equality
$\WH{d} \circ \WH{d} =0$.

We recall cyclic $\AI$-algebras in the gapped filtered case.

\begin{definition}
A filtered gapped $\AI$-algebra $(A,\{m_*\})$ is said to have a {\it cyclic symmetric} inner product if
there exists a skew-symmetric non-degenerate, bilinear map $$<,> : \overline{A}[1] \otimes \overline{A}[1] \to \kk,$$
which is extended linearly over $A$,
such that for all integer $k \geq 0$, $\beta \in G$,
\begin{equation}\label{cyeqnfil}
    <m_{k,\beta}(x_1,\cdots,x_k),x_{k+1}> = (-1)^{K}<m_{k,\beta}(x_2,\cdots,x_{k+1}),x_{1}>.
\end{equation}
 where $K = |x_1|'(|x_2|' + \cdots +|x_{k+1}|')$.
 For short, we will call such an algebra, cyclic (filtered) $\AI$-algebra.
\end{definition}

The discussion in the previous subsection can be extended to the filtered case easily, and
we leave the details to the reader.

\subsection{Hochschild homology of $\AI$-algebra}
We recall the definition of Hochschild and cyclic homology of an $\AI$-algebra
for reader's convenience. We assume that the base ring of $\AI$-algebra contains $\QQ$. 
Let $(A,\{m_k\})$ be a filtered $\AI$-algebra. We denote
$$C^k(A,A) = A[1] \otimes A[1]^{\otimes k}.$$
We will denote its degree $\bullet$ part as $C^k_\bullet(A,A)$.
We define the Hochschild chain complex
\begin{equation}\label{def:hochchain}
C_{\bullet}(A,A) = \HH{\oplus}_{k \geq 0} C^k_\bullet(A,A),
\end{equation}
after completion with respect to energy filtration and 
with the boundary operation $$d^{Hoch} :C_{\bullet}(A,A) \to C_{\bullet+1}(A,A)$$  defined as follows: 
For $v \in A$ and $x_i \in A$, 
$$d^{Hoch} (\UL{v} \otimes x_1 \otimes \cdots \otimes x_k)
= \sum_{\stackrel{0 \leq j \leq k+1 -i}{1 \leq i}} (-1)^{\E_1} \UL{v} \otimes \cdots \otimes x_{i-1} \otimes m_j(x_i,\cdots,x_{i+j-1}) \otimes \cdots \otimes x_k$$
\begin{equation}\label{eq22}
 + \sum_{\stackrel{0 \leq i, j \leq k}{ i+j \leq k}} (-1)^{\E_2} \UL{m_{i+j+1}\big(x_{k-i+1},\cdots,x_k,v, x_1,\cdots,x_{j} \big)} \otimes x_{j+1} \otimes \cdots 
\otimes x_{k-i}
\end{equation}
We underline the module elements and the signs follow Koszul convention:
$$\E_1 =  |v|' + |x_1|' + \cdots + |x_{i-1}|',$$
$$\E_2 = \big(\sum_{s=1}^i |x_{k-i+s}|'\big)\big( |v|' + \sum_{t=1}^j |x_{t}|'\big).$$
Here, we are considering the $\AI$-algebra as $\AI$-module over itself.

Let us recall cyclic homology of $A$.
For the cyclic generator $t_{n+1} \in \ZZ/(n+1)\ZZ$, we define its action on $A^{\otimes (n+1)}$:
$$t_{n+1} \cdot (x_0,x_1,\cdots,x_n) = (-1)^{|x_n|'(|x_0|'+\cdots+|x_{n-1}|')} (x_n,x_0,\cdots,x_{n-1}).$$
Here, we set $t_1$ to be identity on $A$ and write the identity map as $1$. Consider
$N_{n+1}:=1+t_{n+1}+t_{n+1}^2+ \cdots+ t_{n+1}^n$. 

As in the classical case,
we have the natural augmented exact sequence:
$$A^{\otimes (n+1)} \stackrel{1-t_{n+1}}{\longleftarrow} A^{\otimes (n+1)}  \stackrel{N_{n+1}}{\longleftarrow}  A^{\otimes (n+1)}  \stackrel{1-t_{n+1}}{\longleftarrow}  A^{\otimes (n+1)}  \stackrel{N_{n+1}}{\longleftarrow} \cdots .$$   
We consider $\oplus_{n=1}^\infty N_n$ action on $\oplus_{n=1}^\infty A^{\otimes n}$ and denote it as
\begin{equation}\label{symop}
N:C_\bullet(A,A) \mapsto C_\bullet(A,A).
\end{equation}
We can also similarly define $(1-t):C_\bullet(A,A) \mapsto C_\bullet(A,A)$.

The Connes' complex for cyclic homology is defined as
$$(C^\lambda(A,A),d^{Hoch})= (C_\bullet(A,A)/1-t,d^{Hoch}).$$ 
From the invariant-coinvariant relation, one can also define cyclic homology
by $(C_\bullet(A,A))^{cyc}$ considering invariant elements of cyclic action on $C_\bullet(A,A)$
with bar differential $\HH{d}_{bar}$. (see \cite{C5} for more details).

Given an element $(\sum_{\sigma \in \ZZ/ k\ZZ} a_{\sigma(1)}\otimes
\cdots a_{\sigma(k)})$ of
$(C_\bullet(A,A))^{cyc}$, we associate an element $a_1 \otimes \cdots \otimes a_k$ or
$\frac{1}{k}(\sum_{\sigma \in \ZZ/ k\ZZ} a_{\sigma(1)}\otimes \cdots
a_{\sigma(k)})$  of $C^\lambda(A,A)$.
The map may be considered as $N^{-1}$. Here $N^{-1}$ is not
well-defined as a map to $C_\bullet(A,A)$ but is well-defined as a map
to $C^\lambda(A,A)$. 

\section{Generalized counting}
In this section, we show that cyclic symmetry provides a function on
homology theories of A-infinity algebras,
namely, function on Hochschild, cyclic homology of $\AI$-algebra or on
(cyclic) Chevalley-Eilenberg homology of the induced $\LI$-algebra.
(In the preprint of 2006, we called it big cohomology cycles, but they were in fact
well-known homology cycles. See \cite{C5} for more details).

\subsection{Heuristic idea}
We first define $m^+$ operation as in \cite{Fu}.
\begin{definition} For $k \geq 1$, we define
$$m_k^+(x_1,\cdots,x_k) = <m_{k-1}(x_1,\cdots,x_{k-1}),x_k> \;\;\;\in \; \NOV. $$
We may sometimes write $m_k^+$ as $m^+$ since $k$ may be determined by the inputing element.
\end{definition}

Now, we consider the filtered $\AI$-algebra of Lagrangian submanifold defined in \cite{FOOO}.
Consider an intersection problem of a holomorphic disc with $k+1$
marked points intersecting chains $x_1,\cdots,x_{k+1}$
at the corresponding marked points. One observes that with the cyclic
inner product $<,>$ defined from the usual intersection pairing,
the number of such holomorphic discs with prescribed intersection may
be given by $m^+(x_1,\cdots,x_{k+1})$.

\begin{figure}[h]
\begin{center}
\includegraphics[height=1.3in]{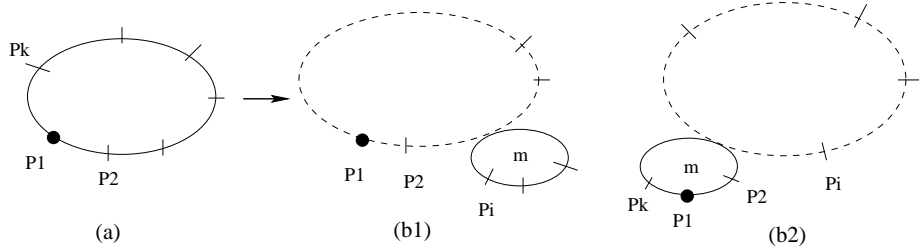}
\caption{Bubbling off and Hochschild boundary operation}
\label{dhat}
\end{center}
\end{figure}

Now, consider a moduli space of holomorphic discs with $k+1$ marked
points intersecting chains $x_1,\cdots,x_{k+1}$ and
its boundary strata from the disc bubbling. We consider the case that the first marked point is specially treated ( in relation to Hochschild theory).  When there exist
a disc bubbling,  there can be two types of bubblings say $(b1)$ or $(b2)$, each of which
corresponds to the intersection problem with 
the expression given by Hochschild differential $d^{Hoch}(x_1,\cdots,x_{k+1})$. 
Here $(b1)$ and $(b2)$ corresponds to the second and third term of (\ref{eq22}). The case
of cyclic homology can be observed by symmetrizing the situation.

Hence the heuristic idea suggests that 
with cyclicity,  intersection with Hochschild boundary $Im (d^{Hoch})$  is realized as a boundary
of a certain moduli space (figure $(a)$), and hence the counting of such intersections (when dimension zero) becomes
zero.  Hence, one observes that  countings do not change after adding $Im (d^{Hoch})$.

But in fact, the argument above
is not entirely correct as it is not yet possible to construct cyclic filtered $\AI$-algebra of singular chains on $L$ directly due to
technical issues of transversality, but Fukaya has constructed a cyclic $\AI$-algebra for the differential forms of $L$.
From now on, we consider general cyclic $\AI$-algebras.

\subsection{Hochschild homology and cyclic symmetry}
We show that $m^+$ is well-defined on $C_\bullet(A,A)/Im(d^{Hoch})$. When restricted to $d^{Hoch}$-cycles, we obtain 
a well-defined function on Hochschild homology of $A$.
\begin{theorem}\label{generalfunction2}
For a cyclic $\AI$-algebra $A$, $m^+ $ gives rise to a well-defined  function on 
 $C_\bullet(A,A)/Im (d^{Hoch})$. In particular $m^+$ defines a function on the Hochschild homology of $A$ and on the Chevalley-Eilenberg homology of the underlying $\LI$-algebra $\WT{A}$
$$m^+:H_\bullet(A,A) \to \NOVO, \;\; \;\; m^+: H^{CE}_\bullet(\WT{A},\WT{A}) \to \NOVO.$$
\end{theorem}
\begin{proof}
For this, we only need to prove that
$$m^+ \big( d^{Hoch}(\UL{x_0} \otimes \cdots \otimes x_{k}) \big) =0$$
From the definition of Hochschild differential (we omit $(-1)^{Kos}$), LHS equals
\begin{equation}\label{eq:g1}
m^+( \sum_{\stackrel{0 \leq j \leq k+1 -i}{1 \leq i}} \UL{x_0} \otimes \cdots \otimes x_{i-1} \otimes m_j(x_i,\cdots,x_{i+j-1}) \otimes \cdots \otimes x_k ) 
\end{equation}
\begin{equation}\label{eq:g2}
 + m^+( \sum_{\stackrel{0 \leq i, j \leq k}{ i+j \leq k}} \UL{m_{i+j+1} \big(x_{k-i+1},\cdots,x_k,x_0, x_1,\cdots,x_{j} \big)} \otimes x_{j+1} \otimes \cdots \otimes x_{k-i} )
\end{equation}
The expression (\ref{eq:g1}) equals (by applying cyclic rotation if necessary)
\begin{equation}\label{eq:g3}
\sum_{\stackrel{0 \leq j \leq k+1-i}{1 \leq i}} <m_{k+1-j}( \cdots,x_0,\cdots,x_{i-1}),
m_j(x_i,\cdots,x_{i+j-1})> 
\end{equation}
This equals the expression in the Proposition \ref{prop:cycplus1}, hence vanishes.

The expression (\ref{eq:g2}) equals 
\begin{equation}\label{eq:g4}
 \sum_{\stackrel{0 \leq i, j \leq k}{i+j \leq k}}
<m_{k-i-j} \big( m_{i+j+1} (x_{k-i+1},\cdots,x_0, \cdots,x_{j} ), x_{j+1},\cdots,x_{k-i-1} \big),x_{k-i}>
\end{equation}
\begin{equation}\label{eq:g5}
=\sum_{\stackrel{0 \leq i, j \leq k}{ i+j \leq k}}
<m_{k - i-j}\big( x_{j+1},\cdots,x_{k-i-1},x_{k-i} \big), m_{i+j+1} \big(x_{k-i+1},\cdots,x_0, \cdots,x_{j} \big)> 
\end{equation}
This also vanishes from the Proposition \ref{prop:cycplus1}.

The statement on Chevalley-Eilenberg homology follows from the Hochschild case:
if we use the notation
$$[x_1,\cdots,x_k]= \sum_{\tau \in S_k} (-1)^{\E(\tau,\vec{x})} x_{\tau(1)}\otimes \cdots \otimes x_{\tau(k)},$$
an element of Chevalley-Eilenberg chain may be written as
$\UL{x_0} \otimes [x_1,\cdots,x_k]$ and its Chevalley-Eilenberg differential also can be understood as
$$d^{CE}(\UL{x_0} \otimes [x_1,\cdots,x_k])
= d^{Hoch}(\UL{x_0} \otimes [x_1,\cdots,x_k]).$$
We refer readers to section 4.1 of \cite{C5} for more details.
\end{proof}

\subsection{The case of cyclic homology}
Note that for cyclic $\AI$-algebra, we have for any $\sigma \in \ZZ/(k+1)\ZZ$,
\begin{equation}\label{eq0}
m^+(x_0,\cdots,x_k) = (-1)^{K(\vec{x})} m^+(x_{\sigma(1)},\cdots,x_{\sigma(k)}).
\end{equation}
From this, we note that $m^+$ vanishes on $Im(1-t)$.
And as cyclic homology is defined by the Connes' complex $(C_\bullet(A,A)/1-t,d_{Hoch})$,
the discussion for Hochschild homology implies also that $m^{+}$ gives
a well-defined function on cyclic homology also.
By symmetrization, one obtains similar statement for cyclic
Chevalley-Eilenberg homology:
\begin{theorem}\label{generalfunction}
For an $\AI$-algebra $A$, $m^+ $ gives rise to well-defined  functions
on the cyclic homology of $A$ and the cyclic Chevalley-Eilenberg homology of the underlying $\LI$-algebra $\WT{A}$
$$m^+:HC_\bullet(A) \to \NOVO, \;\; \;\; m^+: HC^{CE}_\bullet(\WT{A}) \to \NOVO.$$
\end{theorem}

We remark that the $\AI$-equation can be translated into the following equation
\begin{lemma} We have
$$ m^+ \big( \widehat{d} (x_{1} \otimes \cdots  \otimes x_{k}) \otimes x_{k+1} \big) =0.$$
\end{lemma}
\begin{proof}
Recall that $\AI$-equation may be written as $m( \widehat{d}(x_{1} \otimes \cdots  \otimes x_{k}))=0$, and
the lemma follows from the definition of $m^+$.
\end{proof}
\begin{remark}
One might hope to have $$m^+ \big(\widehat{d}(x_{1} \otimes \cdots   \otimes x_{k+1}) \big) =0.$$ but
this does not hold in general. For example, one can easily make an example that
$m^+(\widehat{d}(x_1 \otimes \cdots \otimes x_4)) \neq 0$ with the assumption $m_0=m_1=0$.
\end{remark}

\subsection{Invariance}
We explain the invariance of the function $m^+$ with respect to the cyclic $\AI$-homomorphism.
Let $(A,m_A,<,>_A)$, $(B,m_B,<,>_B)$ be two cyclic $\AI$-algebras and let $h:A \to B$ be an $\AI$-homomorphism.
Then, from $h$ one can define the cohomomorphism $\HH{h}$, and an induced map $\HH{h}:C_\bullet(A,A) \to C_\bullet(B,B)$ which
is a chain map between two Hochschild chain complexes.
\begin{prop}
If $h:A \to B$ is a {\em cyclic} $\AI$-homomorphism, it preserves the value
of $m^{+}$.
Namely, for a Hochschild cycle (or chain) $\alpha$, we have
$$m^+_A(\alpha) =m^+_B(\HH{h}(\alpha))$$
\end{prop}
\begin{proof}
The proof follows from the definition of a cyclic $\AI$-homomorphism (Definition \ref{def:homo}).
Namely
$$m^+_B(\HH{h}(\cdots)) =
<m_B(h(\cdots),h(\cdots),\cdots,h(\cdots)),h(\cdots)>_B$$
$$=<h(\cdots,m_A(\cdots),\cdots),h(\cdots)>_B =
<h_1(m_A(\cdots)),h_1(\cdot)>_B$$
$$=<m_A(\cdots),\cdot>_A$$
Here, first equality follows from the definition of $\HH{h}$, and the second equality is $\AI$-equation on
the LHS, and the third and the last equalities are from the definition of cyclic $\AI$-homomorphism.
\end{proof}

\section{A toric example}
We illustrate the above formalism by considering a specific cyclic cohomology class of the $\AI$-algebra of the Clifford torus $T^2$ in $\CP^2$. 
(This example is inspired by the recent  work of \cite{FOOO2} in particular, their invariant Z, but we do not know the exact correspondence between them.)

One can construct the $\AI$-algebra of the Clifford torus  as follows.
Denote by $A^\bullet$ the differential forms on $L=T^2$ (identified as $\RR^2/\ZZ^2$).
We will consider the Novikov field $\Lambda$ which is obtained by forgetting $e$ parameter from $\NOV$.
Then, the map $m_{k,\beta}$ is defined as follows.
Consider the evaluation map
$$ev=(ev_1,\cdots,ev_k,ev_0):\CM_{k+1}^{main}(L,\beta) \to L^{k+1}.$$ 
For $\rho_1,\cdots,\rho_k  \in A^\bullet$, define
\begin{equation}\label{defit}
m_{k,\beta}(\rho_1,\cdots,\rho_k) = ev_0!(ev_1,\cdots,ev_k)^*(\rho_1 \wedge \cdots \wedge \rho_k).
\end{equation}
In the toric case, $ev_0!$ is well-defined as $ev_0$ is a submersion by  the torus action.

The canonical model $A_{can}$ then satisfies the following conditions:
\begin{enumerate}
\item $A_{can}$ is finite dimensional and $(A_{can},\{m_{can,k}\})$ has a unit $e$.
\item $m_{can,0}(1)$ is a constant multiple of the unit $e$. (Here $e$ is a constant function 1 on $L$).
\item $m_{can,1}=0$.
\end{enumerate}
When we consider cyclic cohomology of $A_{can}$, we replace $m_{can,0}(1)$ by $0$, and consider that of the resulting $\AI$-algebra 
and denote the resulting algebra as $(A, \{m_k\}$ for simplicity. ( We first recall some computations of $m_1$ and $m_2$ operations based on the work of
\cite{C1}, \cite{C2},\cite{CO}, \cite{FOOO1}. We will not mention about the signs and details of the computation, which refer
readers to the above literatures.)

Denote by $e_1=dx, e_2=dy$ the harmonic differential forms on $\RR^2/\ZZ^2 =\{(x,y)\}$ with the flat metric.
Then, the following computation (see Lemma 10.8 of \cite{FOOO1}) only involve holomorphic discs of Maslov index two and has
been computed. Recall that $\beta_0, \beta_1, \beta_2$ are the homotopy classes of 
basic holomorphic discs, and we assume that the symplectic area of $\beta_i$ equals 1 for all $i$.

\begin{equation}
m_{k,\beta_j}(e_{i_1},\cdots,e_{i_k})= \frac{1}{k!} \prod_{s=1}^k (\partial \beta_j \cap e_{i_s})T \\
\end{equation}
In particular, $m_2(e_1,e_1)=m_2(e_2,e_2)=1T$ and $m_2(e_1,e_2) = m_2(e_2,e_1) = (1/2 )T$.
The canonical model  as a vector space over $\Lambda$ is generated by $e_1$ and $e_2$.
Hence, the product on the canonical model defines Clifford algebra
$$e_i \cdot e_j + e_j \cdot e_i = H_{ij} T, \;\; \textrm{where} \;\; H = (\begin{array}{cc} 2 &1 \\ 1&2 \end{array}).$$

The eigenvalues of the Hessian matrix $H$ is 3 and 1, and the we take the following eigenvectors as basis of $A_{can}$.
$$f_1 = \frac{e_1+e_2}{\sqrt{2}},\;\; f_2 = \frac{e_1-e_2}{\sqrt{2}}.$$
Then, it is easy to check that 
$$m_2(f_1,f_1) = \frac{3}{2} T, \;\;m_2(f_2,f_2) = \frac{1}{2} T,$$
$$  m_2(f_1,f_2)=e_1\wedge e_2,\;\;
m_2(f_2,f_1) = -e_1 \wedge e_2.$$
We set $f_{ij} = m_2(f_i,f_j)$ and we have $f_{21} = -f_{12}$.

Recall that $N_3$ is the symmetrization operator for the tensor of length 3. 
 \begin{prop} The following expression $\alpha$
 is a non-trivial cyclic cohomology cycle.
$$\alpha :=  N_3(f_1 \otimes f_1 \otimes f_{12}) + 3 N_3(f_2 \otimes f_2 \otimes f_{12}) + 3N_3(e \otimes f_1 \otimes f_2)T$$
The value of the function $m^+$ equals $18$ which also implies that the above cycle is
non-trivial.
\end{prop}
\begin{remark}
The invariant $Z$ of \cite{FOOO2} in this case is given by $6$, which corresponds to the $m_+$ invariant of $\frac{1}{3}\alpha$.
\end{remark}
\begin{proof}
First, we check that $\alpha$ is indeed a cyclic cohomology cycle.
The unit identity $ m_2(e,x) = (-1)^{deg(x)}m_2(x,e) =x$ will be used several times.

From the $\AI$-formula 
$m \circ \WH{m} (f_2,f_1,f_1)=0$ and $m \circ \WH{m} (f_1,f_2,f_2)=0$,
we obtain that
\begin{eqnarray*}
m_2(f_{12}, f_1) &=& -\frac{3}{2} f_2 T\\
m_2(f_{12}, f_2) &=& \frac{1}{2} f_1 T.
\end{eqnarray*}
Namely, if one expands $m \circ \WH{m} (f_2,f_1,f_1)=0$, 
the non-vanishing terms are
$$ m_2(m_2(f_2,f_1),f_1) + m_2(f_2, m_2(f_1,f_1)) =0$$
and hence
$$m_2(f_{12},f_1) = -m_2(f_{21},f_1) = -(- m_2(f_2,m_2(f_1,f_1))) = m_2(f_2, \frac{3}{2}T\cdot e) = -\frac{3}{2}f_2T
$$

In the same way, one can also check that 
\begin{eqnarray*}
m_2(f_1, f_{12}) &=& -\frac{3}{2} f_2 T\\
m_2(f_2, f_{12}) &=& \frac{1}{2} f_1 T 
\end{eqnarray*}

We use the above computations to show that $\WH{m}_2 \alpha =0$(we omit $\otimes$ and $T$ for simplicity).

First, we have
\begin{equation}\label{e1}
\WH{m}_2 ( e f_1f_2 + f_2ef_1 + f_1f_2e)
= (f_1f_2 - e f_{12}) + ( -f_2f_1 + f_2f_1) + (f_{12}e - f_1f_2)
= f_{12}e - ef_{12}
\end{equation}
Next, we have
$$\WH{m}_2(f_{12}f_1f_1 + f_1f_{12}f_1 + f_1f_1f_{12})$$
$$= N_2(m_2(f_{12},f_1) \otimes f_1) + N_2(m_2(f_1,f_{12})\otimes f_1) 
+N_2(m_2(f_1f_1)\otimes f_{12})$$
Hence,
\begin{equation}\label{e2}
\WH{m}_2(f_{12}f_1f_1 + f_1f_{12}f_1 + f_1f_1f_{12})
= N_2( - \frac{3}{2} f_2\otimes f_1) + N_2(- \frac{3}{2} f_2\otimes f_1) 
+N_2( \frac{3}{2} e\otimes f_{12})
\end{equation}
(Here, as the shifted degree of $f_1$ is zero and shifted degree of $e$ and $f_{12}$ is odd,
we have $N_2(f_1f_2) = f_1f_2+f_2f_1, \; N_2(ef_{12}) = ef_{12} - f_{12}e$).

Similarly, we have
\begin{equation}\label{e3}
\WH{m}_2(f_{12}f_2f_2 + f_2f_{12}f_2 + f_2f_2f_{12})
= N_2( \frac{1}{2} f_1\otimes f_2) + N_2(\frac{1}{2} f_1\otimes f_2) 
+N_2(\frac{1}{2} e\otimes f_{12})
\end{equation}
By taking $ 3 \times (\ref{e1}) + (\ref{e2}) + 3\times (\ref{e3})$, we check that $\WH{m}_2 \alpha =0$.

Now, we show that $m_3 \alpha =0$.
First, note that from the unital property, we have
 $$m_3(ef_1f_2)=m_3(f_2ef_1)=m_3(f_1f_2e)=0.$$
Hence, we only need to compute $m_3(N_3(f_{12}f_1f_1))$ and $m_3(N_3(f_{12}f_2f_2))$.
As before, we use the $\AI$-equation $m \circ \WH{m}(f_2,f_1,f_1,f_1)=0$.
This shows that $$m_3(f_{21},f_1,f_1) = - m_2(m_3(f_2,f_1,f_1), f_1) -m_2(f_2,m_3(f_1,f_1,f_1)).$$
Hence, we have 
$$m_3(f_{12},f_1,f_1) = m_2(\frac{1}{6}\sum_i (\partial \beta_i \cap f_2)(\partial \beta_i \cap f_1)^2 e, f_1) 
+ m_2(f_2, \frac{1}{6}\sum_i (\partial \beta_i \cap f_1)^3e).$$
$$ = \frac{1}{6} \big( ( (\frac{1}{\sqrt{2}})^3 - (\frac{1}{\sqrt{2}})^3)f_1 - ( (-\sqrt{2})^3 + (\frac{1}{\sqrt{2}})^3 + (\frac{1}{\sqrt{2}})^3)f_2 \big)$$
$$=  \frac{\sqrt{2}}{4} f_2. $$
One can check that
$$m_3(f_{12},f_1,f_1) = m_3(f_1,f_{12},f_1) = m_3(f_1,f_1,f_{12}) =  -\frac{\sqrt{2}}{4} f_2.$$

Similarly, we compute $m_3(N_3(f_{12}f_2f_2))$.
$$m_3(f_{12},f_2,f_2) = - m_2(f_1,m_3(f_2,f_2,f_2)) - m_2(m_3(f_1,f_2,f_2),f_2)$$
$$=  \frac{1}{6}(\sum (\partial \beta_i \cap f_2)^3) f_1  - (\sum (\partial \beta_i \cap f_1)(\partial \beta_i \cap f_2)^2 f_2)$$
$$ = -  \frac{1}{6}((\frac{1}{\sqrt{2}})^3 + (\frac{1}{\sqrt{2}})^3)f_2 = -  \frac{\sqrt{2}}{12}f_2$$

One can also check that $m_3(f_{12}f_2f_2) =m_3(f_2f_{12}f_2) = m_3(f_2f_2f_{12})$.
Therefore
$$m_3(\alpha) = \big(  3 \frac{\sqrt{2}}{4} - 3 \times 3 \frac{\sqrt{2}}{12} \big) f_2 = 0.$$
Hence, this finishes the proof that $\alpha$ is a cyclic cohomology cycle.

Now, we compute $m_+$ invariant of $\alpha$.
Consider cyclic inner product satisfying 
$$<f_1,f_2> = - <f_2,f_1> = 1$$
Cyclic property implies that 
$$<f_{12},e> = <m_2(f_1,f_2),e> = <m_2(f_2,e),f_1> = < -f_2,f_1> = 1 = <e,f_{12}>.$$

The following expressions contribute to the $m_+$ invariant of $\alpha$:
$$<m_2(f_1,f_1),f_{12}> + <m_2(f_1,f_{12}),f_1> + <m_2(f_{12},f_1),f_1>$$
$$3<m_2(f_2,f_2),f_{12}> + 3<m_2(f_2,f_{12}),f_2> + 3<m_2(f_{12},f_2),f_2>$$
$$3<m_2(L,f_1),f_2> + 3<m_2(f_2,L),f_1> + 3<m_2(f_1,f_2),L>$$
we can verify that
sum of the a expression equal
$(3/2 +3/2 + 3/2)+ (3/2 + 3/2 +3/2) + (+3 +3 +3) = +18$.
This proves the theorem.
\end{proof}

\bibliographystyle{amsalpha}

\end{document}